\documentclass[12pt]{article}
\usepackage{amsmath,amssymb,latexsym}
\usepackage{longtable}
\usepackage{amsthm}
\usepackage{graphicx}
\usepackage{epsfig}

\newtheorem{theorem}{Theorem}

\newtheorem{proposition}[theorem]{Proposition}

\newtheorem{example}{Example}
\newtheorem{remark}{Remark}

\begin{document}

\begin{center}
{\Large \textbf{On a class of distributions generated by
stochastic mixture of the extreme order statistics of a sample of
size two}}
\end{center}

\begin{center}
{\large \textbf{S. M. Mirhoseini$^{\mathrm{a}}$, Ali
Dolati$^{\mathrm{b}}$ and M. Amini$^{\mathrm{c}}$}}
\end{center}

\begin{center}
\textit{$^{\mathrm{a,c}}$Department of Statistics, Faculty of
Mathematical Science, Ferdowsi University,P.O. Box 91775-1159,
Mashhad, Iran \\ \textrm{mmirhoseini@yazduni.ac.ir}
 \\ $^{\mathrm{b}}$Department of Statistics, Faculty of
Mathematics, Yazd University, Yazd, 89195-741, Iran\\[0pt]\textrm{adolati@yazduni.ac.ir}
\\[0pt]}
\end{center}

\begin{center}
\textbf{Abstract}
\end{center}

\noindent This paper considers a family of distributions
constructed by a stochastic mixture of the order statistics of a
sample of size two. Various properties of the proposed model are
studied. We apply the model to extend the exponential and
symmetric Laplace distributions. An extension to the bivariate
case is considered.
\bigskip

\noindent AMS(2000) Subject classification: 62F15.
\medskip

\noindent \textit{Keywords}: Aging characteristics, copula, hazard
rate function, mixture, stochastic orders; transformation.
\section{Introduction}

Different methods may be used to introduce a new parameter to a
family of distributions to increase flexibility for modeling
purposes. Marshall and Olkin \cite{Marshal97} introduced a method
for adding a parameter to a family of distributions and applied it
to the exponential and Weibull models. Jones \cite{Jones} used the
distribution of order statistics to provide new families of
distributions with extra parameters. The well--known
Farlie--Gumbel-Morgenstern (FGM, for short) family of bivariate
distributions with the given univariate marginal distributions
$F_1$ and $F_2$, is defined by
\begin{equation}
H(x,y)=F_{1}(x)F_{2}(y)\{ 1+\lambda
\bar{F_{1}}(x)\bar{F_{2}}(y)\}, \label{FGM}
\end{equation}
where $\lambda \in \lbrack -1,1]$; see, Drouet-Mari and Kotz
(\cite{Kotz2001}, Chapter 5) for a good review. For a given
univariate cumulative distribution function $F$, the univariate
version of (\ref{FGM}) may be considered as
\begin{eqnarray}
G_{\lambda}[F](x)= F(x)\{1+\lambda \bar{F}(x)\}, \label{m2}
\end{eqnarray}
for all $x$ and $-1\leq \lambda \leq 1$. The family of
distributions defined by (\ref{m2}) is comparable with the
Marshall-Olkin \cite{Marshal97} family of distributions, which
also called the proportional odds model
\cite{Kirmani2001,Marshal97}, given by
\begin{equation}
H(x)=\frac{F(x)}{1-(1-\alpha)\bar{F}(x)}, \quad\ -\infty<x<\infty,
\quad\ \alpha>0.
\end{equation}
Note that $H$ with $0<\alpha<1$ could be written as
\begin{equation}
H(x)=F(x)\sum_{k=0}^{\infty}\{(1-\alpha)\bar{F}(x)\}^k,
\end{equation}
and hence for $0<\lambda<1$, (\ref{m2}) is a first-order
approximation to the proportional odds model.

The aim of the present paper is to investigate different
properties of (\ref{m2}). We first provide a physical
interpretation for this model in Section 3. Some preservation
results of stochastic orderings and aging properties are given in
Section 4. A generalization of the ordinary exponential
distribution which exhibits both increasing and decreasing hazard
rate functions and a skew extension of the symmetric Laplace
distribution are given in Section 5. Bivariate case is discussed
in Section 5. In Section 1 we recall some notions that will be
used in the sequel.

\section{Preliminaries}

Let us recall some notions of stochastic orderings and aging
concepts that will be useful in this paper. Let $X$ be a
continuous random variable with the cdf $F$, the survival function
$\bar{F}=1-F$, the probability density function (pdf) $f$, the
residual life survival function $\bar{F}_t(x)=P(X>x+t|X>t)$ and
the hazard rate function $h_F(x)=f(x)/\bar{F}(x)$. Then $F$ is
said to have: (i) increasing (decreasing) hazard rate IHR (DHR) if
$h_F(x)$ is increasing (decreasing) in $x$; (ii) increasing
(decreasing) hazard rate average IHRA (DHRA) if
$\int_{0}^{t}h_F(x)dx/t$ is increasing (decreasing) in $t$; (iii)
new better (worse) than used NBU (NWU) property if
$\bar{F}_t(x)\leq (\geq) \bar{F}(x)$, for all $x\geq 0$ and
$t\geq0$. The implications
$$
\text{\rm IHR} \Longrightarrow \text{\rm IHRA} \Longrightarrow
\text{\rm NBU} \quad {\rm and} \quad \text{\rm DHR}
\Longrightarrow \text{\rm DHRA} \Longrightarrow \text{\rm NWU},
$$
are well known. See \cite{bar1975} for more detail. The following
definitions will be used for various stochastic comparisons. Let
$F_1$ and $F_2$ be two cdfs with the corresponding pdfs $f_1$ and
$f_2$, the hazard rate functions $h_{F_{1}}$, $h_{F_{2}}$, and the
quantile functions $F^{-1}_1$ and $F^{-1}_2$, respectively, where
$F^{-1}_i=\text{sup}\{x|F_i(x)\leq u\}$, for $0\leq u\leq1$. The
cdf $F_1$ is said to be smaller than $F_2$ in (i) stochastic order
$(F_1\prec_{\rm st}F_2)$ if $F_1(x)\geq F_2(x)$ for all $x$; (ii)
hazard rate order $(F_1\prec_{\rm hr}F_2)$ if $h_{F_{1}}(x)\geq
h_{F_{2}}(x)$ for all $x$; (iii) likelihood ratio order
$(F_1\prec_{\rm lr}F_2)$ if $f_2(x)/f_1(x)$ is non-decreasing in
$x$; (iv) convex transform order $(F_1\prec_{\rm c}F_2)$ if
$F^{-1}_{2}F_1(x)$ is convex in $x$ on the support of $F_1$;  (v)
star order $(F_1\prec_{*}F_2)$ if $F^{-1}_{2}F_1(x)/x$ is
increasing in $x\geq 0$; (vi) superadditive order $(F_1\prec_{\rm
su}F_2)$ if $F^{-1}_{2}F_1(x+y)\geq
F^{-1}_{2}F_1(x)+F^{-1}_{2}F_1(y)$; (vii) dispersive order
$(F_1\prec_{\rm disp}F_2)$ if $F^{-1}_{2}F_1(x)-x$ increases in
$x$. The implications $F_1\prec_{\rm lr}F_2\Longrightarrow
F_1\prec_{\rm hr}F_2\Longrightarrow F_1\prec_{\rm st}F_2$ are well
known. See \cite{Shaked} for an extensive study of these notions.

\section{Genesis of family (2)}
Let $X_{1}$ and $X_{2}$ be two independent and identically
distributed random variables having the survival function
$\bar{F}=1-F.$ For $-1\leq \lambda \leq 1,$ let $Z$ be a Bernulli
random variable, independent of $X_i$s, with
$P(Z=1)=\frac{1+\lambda}{2}$ and $P(Z=0)=\frac{1-\lambda}{2}$.
Consider the stochastic mixture
\begin{equation}
U=ZX_{(1)}+(1-Z)X_{(2)},\label{mrep}
\end{equation}
where $X_{(1)}=\min (X_{1},X_{2})$ and $X_{(2)}=\max
(X_{1},X_{2})$ are the corresponding order statistics of $X_{1}$
and $X_{2}$. Since the distribution functions of $X_{(2)}$ and
$X_{(1)}$ are given by $F_{(2)}(x)=F^{2}(x)$ and
$F_{(1)}(x)=2F(x)-F^{2}(x)$, respectively, then the cdf of $U$,
denoted by $G_{\lambda}[F]$, is given by
\begin{eqnarray}
G_{\lambda}[F](x)&=&\frac{1+\lambda}{2}F_{(1)}(x)+\frac{1-\lambda}{2}F_{(2)}(x)\nonumber\\&=&
F(x)\{1+\lambda \bar{F}(x)\},
\end{eqnarray}
for all $x$ and $-1\leq \lambda \leq 1$. Clearly $G_{0}[F]=F$,
$G_{-1}[F]=F_{(2)}$, and $G_{1}[F]=F_{(1)}$. Since
$G_{\lambda}[F](.)$ is increasing in $\lambda$, we have the
inequality
$$
F_{(2)}(x)\leq G_{\lambda}[F](x)\leq F_{(1)}(x),
$$
for all $x$ and $-1\leq \lambda \leq 1$.

In the following result we show that the transformation (\ref{m2})
is ``unique'', in the sense that given a distribution $F$, this
generates a unique distribution or a family of distributions.

\begin{proposition} Let $F_1$ and $F_2$ be two distribution
functions such that $G_\lambda[F_1]=G_\lambda[F_2]$ for every
$\lambda\in[-1,1]$. Then $F_1=F_2$.
\end{proposition}

\noindent \textbf{Proof.} Suppose that $\lambda> 0$ (the case
$\lambda=0$ is trivial and the case $\lambda<0$ the result could
be proved similar). Then, $G_\lambda[F_1]=G_\lambda[F_2]$, is
equivalent to
\begin{equation}
[F_1(x)-F_2(x)][1-\lambda(F_1(x)+F_2(x)-1)]=0\label{m3},
\end{equation}
for each $x$. Suppose there exist a point $x_0\in R$ such that
---without loss of generality--- $F_1(x_0)<F_2(x_0)$. Then the
equality (\ref{m3}) is equivalent to
$F_1(x_0)+F_2(x_0)=\frac{1}{\lambda}+1$. Since $1\leq
\frac{1}{\lambda}$ and $F_1(x_0)<F_2(x_0)<1$, we must have
$F_2(x_0)>1$. This absurd, so that we conclude that $F_1=F_2$.

\section{Properties}
 The survival function, the probability density function and the hazard rate
 function corresponding to (\ref{m2}) are given by
\begin{equation}
\bar{G}_{\lambda}[F](x)=\bar{F}(x)\{1-\lambda F(x)\},\label{Gbar}
\end{equation}

\begin{equation}
g_{\lambda}[F](x)=f(x)\{1+\lambda(1-2F(x))\}\label{m5}
\end{equation}
and
\begin{equation}
h_{G}(x;\lambda)=\frac{g_{\lambda}[F](x)}{\bar{G}_{\lambda}[F](x)}=h_{F}(x)\left(1+\frac{\lambda
\bar{F}(x)}{1-\lambda F(x)}\right)\label{hr},
\end{equation}
respectively, where, $h_{F}(x)$ is the hazard rate function of
$F$.

It follows from (\ref{hr}) that
\begin{equation*}
\lim_{x\rightarrow -\infty}
h_{G}(x;\lambda)=(1+\lambda)\lim_{x\rightarrow -\infty}h_{F}(x),
\quad\ \lim_{x\rightarrow
\infty}h_{G}(x;\lambda)=\lim_{x\rightarrow \infty} h_{F}(x),
\end{equation*}
\begin{equation*}
h_{F}(x)\leq h_{G}(x;\lambda)\leq (1+\lambda)h_{F}(x), \quad\
(-\infty<x<\infty, 0\leq \lambda\leq 1),
\end{equation*}
\begin{equation*}
\quad\ (1+\lambda)h_{F}(x)\leq h_{G}(x;\lambda)\leq h_{F}(x),
\quad\ (-\infty<x<\infty, -1\leq \lambda\leq 0).
\end{equation*}

Let $\bar{F}_t(x)=\frac{\bar{F}(x+t)}{\bar{F}(t)}$ be the residual
life survival function corresponding to cdf $F$. Then from
(\ref{Gbar}), the residual life survival function of the generated
distribution $G_{\lambda}[F]$, denoted by
$\bar{G}_{\lambda,t}[F](x)$, is given by
\begin{eqnarray}
\bar{G}_{\lambda,t}[F](x)&=&\nonumber\frac{\bar{G}_{\lambda}[F](x+t)}{\bar{G}_{\lambda}[F](t)}
\\\nonumber&=&\bar{F}_t(x)\left(\frac{1-\lambda F(x+t)}{1-\lambda F(t)}\right)\\\nonumber
&=&\bar{F}_t(x)\left(1-\frac{\lambda \bar{F}(t)}{1-\lambda
F(t)}\frac{F(x+t)-F(t)}{\bar{F}(t)}\right)\\\nonumber&=&\bar{F}_{t}(x)\{1-\beta
F_{t}(x)\}\\&=&\bar{G}_{\beta}[F_t](x), \label{rlife}
\end{eqnarray}
where $\beta=\beta(t)=\frac{\lambda \bar{F}(t)}{1-\lambda F(t)}$
and $F_{t}(x)=1-\bar{F}_{t}(x)$. Thus the residual life survival
function of $G_{\lambda}[F]$ is the transformed version of the
residual life survival function of $F$ under (\ref{m2}), with a
new parameter.

By solving the equation
$F(x)\{1+\lambda(1-F(x)\}=G_{\lambda}[F](x)$, with respect to $F$,
one gets
\begin{equation*}
F(x)=\frac{1+\lambda-\sqrt{(1+\lambda)^2-4\lambda
G_{\lambda}[F](x)}}{2\lambda},\label{inv}
\end{equation*}
which gives the the quantile function of $G_{\lambda}[F]$ as
\begin{equation}
G^{-1}_{\lambda}[F](q)=F^{-1}\left( \frac{1+\lambda
-\sqrt{(1+\lambda)^{2}-4\lambda q}}{2\lambda}\right), \quad 0\leq
q\leq1. \quad \label{quantil}
\end{equation}
Note that
$\lim_{\lambda\rightarrow0}G^{-1}_{\lambda}[F](q)=F^{-1}(q)$. In
particular, the median of $G_{\lambda}[F]$ is given by
\begin{equation*}
G^{-1}_{\lambda}[F](0.5)=F^{-1}\left(\frac{1+\lambda
-\sqrt{1+\lambda^{2}}}{2\lambda}\right).
\end{equation*}

\subsection{Stochastic comparisons}

 In this section we provide some results for stochastic orderings and aging
 properties of a given cdf under the transformation (\ref{m2}).

\begin{proposition} For a given cdf $F$, we have

\noindent a) \ \ (i) If $F$ is {\rm IHR} ({\rm IHRA, NBU}) and
$-1\leq \lambda\leq 0$, then $G_{\lambda}[F]$ is {\rm IHR (IHRA,
NBU)}.

\indent (ii) If $F$ is {\rm DHR (DHRA, NWU)} and $0\leq
\lambda\leq 1$, then $G_{\lambda}[F]$ is {\rm DHR (DHRA, NWU)}.

\noindent b) $F\prec_{lr} G_{\lambda}[F]$ (consequently,
$F\prec_{\rm hr} G_{\lambda}[F]$ and $F\prec_{\rm st}
G_{\lambda}[F]$), if $-1\leq\lambda\leq 0$ and
$G_{\lambda}[F]\prec_{\rm lr}F$ (consequently, $F\prec_{\rm hr}
G_{\lambda}[F]$ and $F\prec_{\rm st} G_{\lambda}[F]$) if
$0\leq\lambda\leq1$.

\noindent c)the parametric family $\{G_{\lambda}[F]\}$ of
distributions is decreasing in $\lambda$ in the likelihood ratio
order. Consequently, $G_{\lambda}[F]$ is decreasing in the hazard
rate and stochastic orders.
\end{proposition}

\begin{proposition} Suppose that $F_1$ and $F_2$
be two given CDFs such that $F_1\prec_{\rm st} F_2$. Then
$G_{\lambda}[F_1]\prec_{\rm st} G_{\lambda}[F_2]$ for every
$\lambda\in [-1,1]$.
\end{proposition}

\noindent \textbf{Proof.} Since $F_1\prec_{\rm st}F_2$ implies
that $F_1(x)\geq F_2(x)$ and $\bar{F}_1(x)\leq \bar{F}_2(x)$, for
all $x$; we have
$G_{\lambda}[F_1](x)=F_1(x)\{1+\lambda\bar{F}_1(x)\}\geq
F_2(x)\{1+\lambda\bar{F}_2(x)\}=G_{\lambda}[F_2](x)$ for $\lambda<
0$ and $\bar{G}_{\lambda}[F_1](x)=\bar{F}_1(x)\{1-\lambda
\bar{F}_1(x)\}\leq \bar{F}_2(x)\{1-\lambda
F_2(x)\}=\bar{G}_{\lambda}[F_2](x)$ for $\lambda
>0$, which completes the proof.

\begin{proposition} Let $F_1$ and $F_2$ be two given
cdfs and let $G_{\lambda}[F_1]$ and $G_{\lambda}[F_2]$ be their
transformed versions using (\ref{m2}). Then

$$
F_1\prec_{\rm order}F_2\Rightarrow G_{\lambda}[F_1]\prec_{\rm
order}G_{\lambda}[F_2],
$$

where $\prec_{\rm order}$ is any one of the orders $\prec_{\rm
c}$,  $\prec_{*}$ $\prec_{\rm su}$ and $\prec_{\rm disp}$.
\end{proposition}

\noindent \textbf{Proof.} From (\ref{quantil}) it is easy to see
that
\begin{equation*}
G^{-1}_{\lambda}[F_2]\left(G_{\lambda}[F_1](x)\right)=F^{-1}_2\left(F^{-1}_1(x)\right),
\end{equation*}
for all $x$, which gives the required result.

\subsection{A symmetry property}
The transformation map (\ref{m2}) can be applied to any symmetric
or asymmetric distribution. The following result shows the effect
of this transformation on the symmetry property of the parent
distribution.

\begin{proposition} Let $X$ with the cdf $F$, be a symmetric
random variable about zero (i.e., $X$ and $-X$ have the same
distribution) and let $Y_{\lambda}$ be a random variable
distributed according to $G_{\lambda}[F]$, $-1\leq\lambda\leq 1$.
Then $Y_{-\lambda}$ and $-Y_{\lambda}$ have the same distribution.
\end{proposition}

\noindent \textbf{Proof.} Since $X$ is symmetric about zero, then
$F(x)=1-F(-x)=\bar{F}(-x)$ for all $x$. From (\ref{m2}) and
(\ref{Gbar}) we have
\begin{eqnarray*}
P(-Y_{\lambda}\leq y)&=&\nonumber\bar{G}_{\lambda}[F](-y)\\
&=&\nonumber\bar{F}(-y)\{1-\lambda F(-y)\}\\&=&\nonumber
F(y)\{1-\lambda
\bar{F}(y)\}\\&=&\nonumber G_{-\lambda}[F](y)\\
&=&\nonumber P(Y_{-\lambda}\leq y),
\end{eqnarray*}
which completes the proof.

\section{Examples}
\subsection{The transformed exponential distribution}

 In particular case that $F$ is an exponential
distribution with the parameter $\theta$, the two--parameter
distribution generated using (\ref{m2}) has the cdf
\begin{equation}
G(x;\lambda,\theta)=(1-e^{-\theta x})(1+\lambda e^{-\theta x}),
\quad\ x,\theta>0, -1\leq \lambda\leq1, \label{ME}
\end{equation}
and the corresponding density function
\begin{equation}
g(x;\lambda,\theta)=\theta e^{-\theta x}\{1+\lambda (2e^{-\theta
x}-1)\}.
\end{equation}

For the density function $g$, we have that $\log g(x;\lambda
,\theta )$, is concave for $-1\leq \lambda \leq 0$ and convex for
$0\leq \lambda \leq 1.$ As a result for $0\leq \lambda \leq 1$,
$g(x;\lambda,\theta)$ is decreasing, and for $-1\leq \lambda < 0$,
$g(x;\lambda,\theta)$ is unimodal. By solving the equation $d\log
g(x;\lambda ,\theta )/dx=0$, it is readily verified that the
density function $g(x;\lambda ,\theta )$ has the mode equal to
zero for $\lambda
>-\frac{1}{3}$ and $-\frac{1}{\theta }\ln ( \frac{\lambda-1
}{4\lambda})$ for $\lambda<-\frac{1}{3}$.

From (\ref{hr}), the hazard rate function of this distribution is
given by
\begin{equation*}
h(x;\lambda ,\theta )=\frac{\theta \{ 1+\lambda (2e^{-\theta
x}-1)\} }{1+\lambda (e^{-\theta x}-1)}.
\end{equation*}
It may be noticed that while the exponential distribution has a
constant hazard rate function, the generated cdf $G$, has
increasing hazard rate for $-1\leq\lambda < 0,$ and decreasing
hazard rate for $0< \lambda \leq1$, which follows using the
log-convexity and the log-concavity of the density function.

From (\ref{rlife}) the residual life survival function
corresponding to (\ref{ME}), is given by
\begin{equation}
\bar{G_{t}}(x;\lambda ,\theta )=e^{-\theta x}\{ 1+\beta
(e^{-\theta x}-1)\}, \label{rlifeME}
\end{equation}
where $\beta=\beta(t)=\lambda e^{-\theta t}\{ 1+\lambda
(e^{-\theta t}-1)\}^{-1}$. The limit distribution as
$t\rightarrow\infty$ is an ordinary exponential distribution
because the limit of $\beta(t)$ is 0.

From (\ref{rlifeME}) the mean residual life function of a random
variable $X$ having cdf (\ref{ME}), could be obtained as
\begin{eqnarray*}
m(t;\lambda ,\theta ) &=&\nonumber E(X-t|X>t) \\
&=&\nonumber\int_{0}^{\infty}\bar{G_{t}}(x;\lambda ,\theta )dx\\
&=&\frac{1+\lambda (\frac{1}{2}e^{-\theta
t}-1)}{\theta[1+\lambda(e^{-\theta t}-1)]},
\end{eqnarray*}
which is increasing in $t$ for $0\leq\lambda\leq 1$ and decreasing
for $-1\leq\lambda\leq 0$, with $\lim_{t \rightarrow \infty
}m(t;\lambda ,\theta )=1/\theta=E(X;0,\theta)$ and $
\lim_{t\rightarrow 0}m(t;\lambda ,\theta
)=(2-\lambda)/2\theta=E(X;\lambda,\theta)$; and hence
$$
\frac{1}{\theta}  \leq m(t;\lambda ,\theta )\leq
\frac{2-\lambda}{2\theta } \text{ \ \ \ } (-1\leq \lambda \leq 0),
$$
and
$$
\text{ \ \ \  } \frac{2-\lambda}{2\theta }\leq m(t;\lambda ,\theta )\leq\frac{1}{\theta} \text{ \  \ }%
(0\leq \lambda \leq 1).
$$

The moment generating function of this distribution is given by
\begin{equation*}
M(t)=E(e^{tX})=\frac{\theta \{2\theta -(1+\lambda )t\}}{(\theta
-t)(2\theta -t)}.
\end{equation*}
By straightforward integration the raw moments are found to be
\begin{equation*}
E(X^{r})=\frac{(1+\lambda(2^{-r}-1))r!}{\theta ^{r}},
\end{equation*}
for $r\in N.$

Since for the exponential distribution with the parameter $\theta$
we have $F^{-1}(q)=-\frac{1}{\theta}\text{ln}(1-q)$, $0<q<1$, then
from (\ref{quantil}) the quantile function of the generated
distribution is given by
\begin{equation*}
G^{-1}(q)=-\frac{1}{\theta }\ln \left( \frac{\lambda
-1+\sqrt{(1+\lambda)^{2}-4\lambda q}}{2\lambda }\right).
\end{equation*}
Note that if $\lambda \rightarrow 0,$ then $G^{-1}(q)\rightarrow
-\frac{1}{\theta }\ln (1-q).$

It may be noticed that for the generated distribution,
median($X$), mode($X$) and $E(X)$ are all decreasing in $\lambda$,
 $\theta$ and $\text{mod}(X)\leq \text{median}(X)\leq E(X)$.

\subsection{A class of skew--Laplace distributions}

The classical symmetric Laplace distribution has the pdf
\begin{equation*}
f(x;\theta)=\frac{1}{2\theta}e^{-\frac{|x|}{\theta}},
\end{equation*}
and cdf
\begin{equation*}
F(x;\theta)=\left\{
\begin{array}{ll}
\frac{1}{2}e^{\frac{x}{\theta}}, & \mbox{{\rm} $x\leq 0$}, \\
\noalign{\smallskip}1-\frac{1}{2}e^{\frac{-x}{\theta}},&
\mbox{{\rm } $x\geq 0,$}
\end{array}
\right.\label{SLD}
\end{equation*}
where $-\infty<x<\infty$ and $\theta>0$. The symmetric Laplace
distribution has been used as an alternative to the normal
distribution for modeling heavy tails data. Different forms of the
skewed Laplace distributions have been introduced and studied by
various authors. Recently, Kozubowski and Nadarajah
\cite{Kozub2008} identified over sixteen variations of the Laplace
distribution. In the following we propose a new version of the
skewed Laplace distribution using (\ref{m2}). The cdf and pdf of
the generated model are given by
\begin{equation*}
G_{\lambda}(x;\theta)=\left\{
\begin{array}{ll}
\frac{1}{2}e^{\frac{x}{\theta}}\{1+\lambda(1-\frac{1}{2}e^{\frac{x}{\theta}})\}, & \mbox{{\rm} $x\leq 0$}, \\
\noalign{\smallskip}1-\frac{1}{2}e^{\frac{-x}{\theta}}\{1-\lambda(1-\frac{1}{2}e^{\frac{-x}{\theta}})\},&
\mbox{{\rm } $x\geq 0,$}
\end{array}
\right.\label{SLG}
\end{equation*}
and
\begin{equation*}
g_{\lambda}(x;\theta)=\left\{
\begin{array}{ll}
\frac{1}{2\theta}e^{\frac{x}{\theta}}\{1+\lambda(1-e^{\frac{x}{\theta}})\}, & \mbox{{\rm} $x\leq 0$}, \\
\noalign{\smallskip}\frac{1}{2\theta}e^{\frac{-x}{\theta}}\{1-\lambda(1-e^{\frac{-x}{\theta}})\},&
\mbox{{\rm } $x\geq 0,$}
\end{array}
\right.\label{SLg}
\end{equation*}
respectively. The moment generating function of $G_{\lambda}$, is
given by
\begin{equation*}
M(t)=\frac{1-\lambda \theta t}{1-(\theta t)^2}+\frac{\lambda
\theta t}{4-(\theta t)^2},
\end{equation*}
and the raw moments are found to be
\begin{equation*}
E(X^r)=\left\{
\begin{array}{llll}
\frac{r!\lambda\theta^{r}(1-2^{r+1})}{2^{r+1}}, & \mbox{{\rm if $r$ is odd, }}\\
\noalign{\smallskip}r!\theta^{r},& \mbox{{\rm if $r$ is even. }}
\end{array}
\right.\label{Mom}
\end{equation*}
The expectation, variance, skewness and the kurtosis are given by
\begin{eqnarray*}
\text{E}(X)&=&-\frac{3}{4}\lambda\theta , \\
\text{Var}(X)&=&\theta^{2}(1-\frac{9}{16}\lambda^{2}), \\
\text{Skewness}(X)&=&\frac{18\lambda(4+3\lambda^{2})}{(9\lambda^{2}-32)\sqrt{32-9\lambda^{2}}}, \\
\text{Kurtosis}(X)&=&\frac{6144-243\lambda^{4}-2592\lambda^{2}}{(32-9\lambda^{2})^{2}}.
\end{eqnarray*}
It may be noticed that the skewness of $G_{\lambda}$ is decreasing
in $\lambda$, and then $-1.1423\leq \text{Skewness}(X)\leq
1.1423$. It is positive for $-1\leq\lambda\leq 0$, and negative
for $0\leq\lambda\leq 1$.

\section{Bivariate case}
\subsection{Construction}

 A large number of bivariate distributions
have been proposed in literature. A very wide survey on bivariate
distributions are given in \cite{Bala2009} and \cite{Kotz2000}.
The method used to construct the family of distributions given by
(\ref{m2}) also lends itself well to the construction of bivariate
distributions whose univariate marginal cdf are of the form
(\ref{m2}).

\begin{proposition} Let $F$ be a bivariate cdf with
the univariate
 marginal cdfs $F_1$, $F_2$ and the associated survival function $\overline{F}(x,y)=1-F_1(x)-F_2(y)+F(x,y)$.
 Then, for every $-1\leq\lambda\leq1$, the
 function $G_{\lambda}:\text{R}^2\rightarrow[0,1]$, defined by
\begin{equation}
G_{\lambda}(x,y)=(1+\lambda)\left(F_1(x)F_2(y)+F(x,y)\overline{F}(x,y)\right)-\lambda
F^{2}(x,y), \label{Bivariate}
\end{equation}
is a bivariate cdf with the univariate marginal distributions
\begin{equation}
G_1(x)=F_1(x)\{1+\lambda \overline{F}_1(x)\} \quad \text{and}
\quad G_2(y)=F_2(y)\{1+\lambda \overline{F}_2(y)\}. \label{uni}
\end{equation}
\end{proposition}

\noindent \textbf{Proof.} To prove this, let $(X_1,Y_1)$ and
$(X_2,Y_2)$ be two independent random vector having common
bivariate cdf $F$ and the univariate
 marginal cdfs $F_1$ (of $X_i$) and $F_2$ (of $Y_i$), $i=1,2$.
Let $X_{(1)},X_{(2)}$ and $Y_{(1)},Y_{(2)}$ be their corresponding
order statistics. For $-1\leq\lambda\leq1$, consider the random
pair $(V_1,V_2)=(X_{(1)},Y_{(1)})$ with probability
$\frac{1+\lambda}{2}$ and $(V_1,V_2)=(X_{(2)},Y_{(2)})$ with
probability $\frac{1-\lambda}{2}$. Then, it is straightforward to
verify that $(V_1,V_2)$ have the joint cdf (\ref{Bivariate}) with
$G_\lambda(x,\infty)=G_1(x)$ and $G_\lambda(\infty,y)=G_1(y)$.

Note that the special case $F(x,y)=F_1(x)F_2(y)$, the cdf
(\ref{Bivariate}) reduces to
\begin{equation*}
G_{\lambda}(x,y)=F_1(x)F_2(y)\left\{F_1(x)F_2(y)+(1+\lambda)(\overline{F}_1(x)+\overline{F}_2(y))\right\},\label{FGM2}
\end{equation*}
which may serve as a competitor to the FGM family of distributions
with the univariate margins of the form (\ref{uni}).

\subsection{Underlying copula}
A bivariate distribution $F$ can be written in the form $F(x, y)
=C\{F_1(x),F_2(y)\}$, where $C$ is the copula associated with $F$;
see Nelsen \cite{Nelsen2006} for more detail. The function
$\hat{C}$ defined by
$\hat{C}(u,v)=u+v-1+C(1-u,1-v)=\overline{C}(1-u,1-v)$, is the
survival copula associated with $C$ and, moreover
$\overline{F}(x,y)=\hat{C}\{\overline{F}_1(x),\overline{F}_2(y)\}=\overline{C}\{F_1(x),F_2(y)\}$.

The following result shows the relationship between the copula
associated with the baseline cdf $F$ and the copula of generated
cdf $G_{\lambda}$.

\begin{proposition} Let $C_{\lambda}$ be the copula
of the cdf $G_{\lambda}$ defined by (\ref{Bivariate}) and let $D$
be the copula of the baseline cdf $F$. Then
\begin{eqnarray}
C_{\lambda}(\psi(u),\psi(v))&=&(1+\lambda)\left\{uv+D(u,v)\overline{D}(u,v)\right\}-\lambda
D^2(u,v), \label{cop}
\end{eqnarray}
for all $0<u,v<1$, and $-1\leq\lambda\leq1$, where
$$
\psi(t)=t+\lambda t(1-t), \quad 0<t<1.
$$
\end{proposition}

\noindent \textbf{Proof.} Notice that using the definition of a
copula, the bivariate cdf ({\ref{Bivariate}) can be rewritten as
\begin{eqnarray}
C_{\lambda}\{G_1(x),G_2(y)\}&=&(1+\lambda)\left(F_1(x)F_2(y)+D\{F_1(x),F_2(y)\}\overline{D}\{F_1(x),F_2(y)\}\right)\nonumber\\&-&\lambda
D^{2}\{F_1(x),F_2(y)\}, \label{cop2}
\end{eqnarray}
where $G_1$ and $G_2$ are given by (\ref{uni}). By applying the
transformations $u=F_1(x)$ and $v=F_2(y)$ on both sides of
(\ref{cop2}) we readily obtain the required result.

\begin{remark} Note that if the baseline copula $D$ is symmetric, i.e., $D(u,v)=D(v,u)$, for all $u,v\in (0,1)$, then the generated copula $C_{\lambda}$
defined in (\ref{cop}) is symmetric.
\end{remark}

\begin{proposition} The family of copulas $\{C_{\lambda}\}$ defined in (\ref{cop}) is positively ordered for all $-1<\lambda\leq 1$
and any baseline copula $D$; i.e., $C_{\lambda_{1}}(u,v)\geq
C_{\lambda_{2}}(u,v)$ for all $u,v\in(0,1)$ whenever $\lambda
_{1}\geq \lambda_{2}.$
\end{proposition}

\noindent \textbf{Proof.} For any two
 copulas $C_{\lambda_1}$ and $C_{\lambda_2}$ of the form (\ref{cop}), one has
$$C_{\lambda_1}\{\psi_{\lambda_1}(u),\psi_{\lambda_1}(v)\}-C_{\lambda_2}\{\psi_{\lambda_2}(u),\psi_{\lambda_2}(v)\}=(\lambda_1-\lambda_2)\{D(u,v)(1-u)(1-v)+uv(1-D(u,v))\},$$
where $\psi_{\lambda}(t)=t+\lambda t(1-t)$, $0<t<1$. Since
$D(u,v)(1-u)(1-v)+uv(1-D(u,v))\geq 0$ for all $u,v\in(0,1)$ and
the function $\psi_{\lambda}(t)=t+\lambda t(1-t)$ is increasing in
$\lambda$ for all $t\in(0,1)$, it is easy to see that
$$C_{\lambda_1}(u,v)-C_{\lambda_2}(u,v)\geq 0,$$ for all $u,v\in
(0,1)$ and $\lambda_1\geq\lambda_2$, which completes the proof.
\bigskip

 We now consider some special cases.

\begin{example}
For the special case that $D(u,v)=uv$, i.e.,
$F(x,y)=F_1(x)F_2(y)$, we have
\begin{equation*}
C_{\lambda}\{\psi(u),\psi(v)\}=uv\{uv+(1+\lambda)(2-u-v)\}.
\end{equation*}
\end{example}

\begin{example}
Suppose that $D=M$, where $M(u,v)={\rm min}(u,v)$, is the
Fr\'{e}chet--Hoeffding upper bound copula (see \cite{Nelsen2006});
which means that the baseline cdf $F$ is the cdf of two perfect
positive dependent random variable $X$ and $Y$. Since
$\overline{M}(u,v)=M(1-u,1-v)$ for every $u,v\in(0,1)$, it is easy
to verify that $M(u,v)\overline{M}(u,v)=M(u,v)-uv$. By applying
(\ref{cop}) to $M$, from the fact that for non--decreasing
function $\psi$,
$\text{min}\{\psi(u),\psi(v)\}=\psi(\text{min}(u,v))$ we obtain
\begin{eqnarray*}
C_{\lambda}(\psi(u),\psi(v))&=&M(u,v)\{1+\lambda(1-M(u,v)\}\nonumber\\&=&\psi\{M(u,v)\}\nonumber\\&=&M(\psi(u),\psi(v)),
\end{eqnarray*}
that is for all $\lambda \in (-1,1)$,
$$
C_{\lambda}(u,v)=M(u,v).
$$
Thus the functional transformation (\ref{Bivariate}) preserves the
perfect dependence of the parent distribution.
\end{example}

\begin{example}
Suppose that $D=W$, where $W(u,v)={\rm max}(u+v-1,0)$, is the
Fr\'{e}chet--Hoeffding lower bound copula (see \cite{Nelsen2006});
which means that the baseline cdf $F$ is the cdf of two perfect
negative dependent random variable $X$ and $Y$. It is easy to
verify that $\overline{W}(u,v)W(u,v)=0$, for every $u,v\in(0,1)$.
Thus (\ref{cop}) gives
$$
C_{\lambda}(\psi(u),\psi(v))=uv+\lambda\{uv-W^2(u,v)\}.
$$
\end{example}

%\section{Concluding and Remarks}
%Several parametric family of univariate/ bivariate distributions have been proposed in the literatures, see Balakrishnan and Lai(2009)\cite{Bala2009}.In this paper, we have proposed a parametric family of univariate distributions constructed by a stochastic mixture of the Extreme Order Statistics of a sample of size two.\\
%We developed its properties including some results for stochastic comparisons, symmetry, skewness in special cases.\\
%Moreover, the method used to construct the family of distribution function given in (\ref{m2}) has been applied to the construction of bivariate distribution functions with given univariate marginal distribution function of the form (\ref{m2}).\\
%Some properties of the family of bivariate distributions given in (\ref{Bivariate}), such as, copula function; concordance ordering are obtained.\\
%The results for the dependence structure, Stochastic comparison in bivariate case, multivariate form, and fitting this family to data sets of the family cdf introduced in (\ref{Bivariate}) are our work at near future.

\section{Discussion}
We have introduced a method for constructing a new family of
distributions from any given one. We deliberately restricted our
attention to the study of some general properties of the proposed
model in univariate as well as the bivariate case. The attentive
reader will agree that the construction presented here leaves room
for more studies beyond what accomplished in this work. In our
next investigation we aim to make deeper contributions to the
distribution theory connected to the bivariate case.
\bigskip

\noindent \textbf{\large Acknowledgements} \smallskip

The authors would like to thank an anonymous referee for his/her
valuable suggestions on an earlier version of this paper.

\end{document}